\documentclass[10pt,twoside]{article}
\usepackage{amsmath}
\usepackage{amsfonts}
\usepackage{ecphead}

\usepackage{setspace} 
\def\proof{\noindent{\bf Proof:}\hskip10pt}        
\def\QED{\hfill $\Box$}

\font\tenmath=msbm10 scaled 1200
\font\sevenmath=msbm7 scaled 1200
\font\Fivemath=msbm5 scaled 1200
\newfam\mathfam \textfont\mathfam=\tenmath
\scriptfont\mathfam=\sevenmath \scriptscriptfont\mathfam=\Fivemath

\def\R{\mathbb{R}}
\def \1{1 \mkern -6mu 1} 
\def\N{\mathbb{N}}

\def\P{\mathbb{P}}
\def\F{{\bf F}} 
\def\f{\mathcal{F}}

\def\R{\mathbb{R}}

\def \n{^{(n)}}

\def \e{{\rm e}}

\def \X{{\bf X}}

\def \Y{{\bf Y}}
\def \d{{\rm d}}

\def \t{{\mathcal T}}
\def \r{{\mathcal R}}

\begin{document}
\input{paper-head.sty}

\newtheorem{theorem}{Theorem}
\newtheorem{definition}{Definition}
\newtheorem{proposition}{Proposition}
\newtheorem{lemma}{Lemma}
\newtheorem{corollary}{Corollary}
  
\begin{section}{Introduction}

Aging refers to systems that seem to reach a statistical equilibrium in a certain regime depending on two time-scales: the age of the system and the duration of the observation. If the observation scale is much shorter that the age of the system scale,
then the effect is, more precisely, called subaging.
It arises in a variety of models in random media; see for instance Ben Arous and \v{C}ern\'y \cite{BC} and references therein. 
The purpose of this note is to point out that subaging occurs
in a process which evolves by fragmentation and coagulation, and is related to the work of Pitman \cite{Pi} on coalescent random forests. 
It has its source in an elementary random urn dynamic which we now describe.

\subsection{A two-time-scale phenomenon}

Consider two urns, say $A$ and $B$, and assume that at the initial time $A$ contains $n$ balls and $B$ is empty. At each step $k=1, 2, \ldots$, we flip a fair coin.
If head comes up then we pick a ball in $A$ uniformly at random and place it in $B$ (provided of course that $A$ is not empty, else we do nothing). Similarly, if tail comes up then we pick a ball in $B$ uniformly at random and place it in $A$ (provided of course that $B$ is not empty, else we do nothing). Fix $s,t>0$ arbitrarily.
By the invariance principle, the numbers of balls in $B$ after $\lfloor tn \rfloor$ steps and after $\lfloor tn+ s\sqrt n \rfloor$ steps are both 
close to $\sqrt n R_t$, where $R_t$ is a reflecting Brownian motion evaluated at time $t$, i.e. has the distribution
of the absolute value of an ${\mathcal N}(0,t)$-variable. Precisely because $B$ contains about $\sqrt n$ balls in that period, each ball in $B$ after $\lfloor tn \rfloor$ steps has a probability bounded away from $0$ and $1$ to be selected and placed back in $A$ during the next $\lfloor s\sqrt n \rfloor$ steps. In other words, even though the number of balls in $B$ remains essentially unchanged between $\lfloor tn \rfloor$ and $\lfloor tn+s\sqrt n \rfloor$ steps, the precise contain of urn $B$ evolves significantly.

Now imagine a stochastic process governed by the contain of urn $B$. Suppose that the one-dimension distributions of this process can be renormalized as $n\to\infty$ in such a way that they have a non-degenerate limit, say $\mu_r$, when there are approximately $r\sqrt n$ balls in $B$.   
Let $n$ be large, $t>0$ fixed and  let $s$ vary. The process observed after $\lfloor tn+s\sqrt n \rfloor$ steps then seems to be in statistical equilibrium, in the sense that  its one-dimensional distributions do not change much when $s$ increases. More precisely the almost equilibrium law 
can be expressed as the mixture
$$\int_{0}^{\infty}\mu_r\P(R_t\in\d r)\,.$$
However this is only a pseudo-stationarity  
as this almost equilibrium  distribution depends on the parameter $t$.

The rest of this work is devoted to the rigorous analysis of this two-time-scale phenomenon in the special case when the stochastic process alluded above 
is a fragmentation-coagulation process induced by a natural modification of Pitman's coalescing random forests \cite{Pi}.

\subsection{Fragmentation-coagulation of random forests}
Pitman considered the following random dynamics: one first picks a tree with $n$ labeled vertices uniformly at random, and one deletes its edges one after the other, also uniformly at random. At each step the tree-component containing the edge which is deleted splits into two smaller trees, so after $k\leq n-1$ steps, we obtain a forest with $k+1$ trees, and 
the process terminates after $n-1$ steps when all vertices are isolated.
 Pitman's motivation for considering these edge-deletion  dynamics stems from its connexion with 
the additive coalescence. The latter governs the evolution of a particle system in continuous times in which pairs of particles coalesce with a rate proportional to the sum of their masses. More precisely, Pitman pointed out that the process of the sizes of the tree-components in the forests resulting from the edge-deletion dynamics is a Markovian fragmentation chain, and that 
time-reversal yields the discrete-time skeleton of the additive coalescent started from $n$ monomers (i.e. $n$ atoms each having a unit mass).

On the one hand, it is well-known from the work of Aldous \cite{Al} that 
the uniform random tree on a set of $n$ vertices can be rescaled (specifically edges by a factor $1/\sqrt n$ and masses of vertices by a factor $1/n$)  and then converges weakly as $n\to\infty$ towards the Continuum Random Tree (for short, CRT). 
On the other hand, Evans and Pitman \cite{EP} proved that
 the additive coalescent  started from a large number of monomers
possesses a scaling limit, which is known as the standard additive coalescent.
Roughly speaking, Aldous and Pitman  \cite{AP} put the pieces together: they showed that the properly rescaled edge-deletion procedure on finite trees converges weakly to a Poissonian logging of the CRT. The latter induces the CRT fragmentation of masses, denoted here by $(\F_t)_{t\geq 0}$, and in turn this yields the standard additive coalescent upon time-reversal.

We consider in this work an evolution which combines  edge-removal and edge-replacement, and should be viewed as an avatar of the urn dynamic depicted in the first part. 
In this direction it will be convenient to use marks on edges, agreeing that
a mark on an edge means that this edge has been removed, while the absence of mark means that this edge is present. For each fixed $n$, we first pick a tree on $n$ vertices uniformly at random and declare that initially all edges have no mark. At each step we flip a fair coin. If head comes up then we put a mark on one edge chosen uniformly at random amongst the un-marked edges (provided of course that there still remains at least one edge without mark; else we do nothing), while if tail comes up then we erase the mark 
of one edge chosen uniformly at random amongst the marked edges  (provided of course that there exists such an edge; else we do nothing). 
For every integer $k\geq 0$, we denote by $\f\n_k$ the random forest which results
from removing the marked edges after the $k$-th step  
and by $\X\n_k$ the collection of the sizes of the tree-components in $\f\n_k$ rescaled by a factor 
$1/n$ and ranked in the decreasing order. In a technical jargon,  $\X\n_k$ is a random mass-partition, that is a decreasing sequence  of positive real numbers with sum $1$. Plainly the addition of a mark corresponds to a fragmentation event and erasure  to a coalescence; in other words the chain $(\X\n_k)_{k\in\N}$ evolves by fragmentation and coagulation.

Our purpose is to investigate the asymptotic behavior in distribution of $\X\n_k$ as $k,n\to \infty$. This is easy as far as only one-dimensional distributions are concerned. Indeed Donsker's invariance principle implies that when $n$ is large,  the number of marks after $\lfloor t n\rfloor$ steps is about $\sqrt n R_t$, where $(R_t: t\geq 0)$ is a reflected Brownian motion. It then follows from results  in Aldous and Pitman \cite{AP} that for every fixed $t>0$, 
$ \X\n_{ \lfloor t n\rfloor} $ converges weakly as $n\to\infty$ towards $\F_{R_t}$,  the CRT mass-fragmentation observed at the independent random time $R_t$.

 It turns out that things are more subtle for finite-dimensional marginals.
 The simple rescaling of times by a factor $n$ is too crude and the asymptotic behavior in law is better revealed in the finer regime $tn+s\sqrt n$. Indeed our main result says that
for every fixed
$t>0$,  the process $( \X\n_{ \lfloor t n + s \sqrt n\rfloor})_{s\in\R}$
converges weakly in the sense of finite-dimensional marginals
towards some non-degenerate stationary process with stationary law given by
the distribution of $\F_{R_t}$.
Thus the fragmentation-coagulation process $\X\n$ exhibits  subaging, in the sense that this process seems to reach a statistical equilibrium in the regime $tn+s\sqrt n$
when $n$ is large and $t$ fixed.

The plan of the rest of this note is as follows. Our main result is stated and proved in 
Section 2 after recalling some notions on the CRT. Finally Section 3 is devoted to some comments, complements and open questions.
Our approach owes much to the construction by Aldous and Pitman of the standard additive coalescent via Poissonian cuts on the  skeleton of a Continuum Random Tree.

\end{section}

\begin{section}{Main result}

We start by recalling some elements on the CRT, refering the reader to Aldous \cite{Al}, Evans \cite{Ev} and Le Gall \cite{LG} for background, and then state our main result.

Let $\t$ be a Brownian CRT; in particular $\t$ is almost surely a compact metric space which has the structure of a real tree. Extremities of $\t$ are called leaves; in other words
$a\in\t$ is a leaf means that if $a$ lies on some path joining two points $b,c$ in $\t$, then necessarily $a=b$ or $a=c$. The subset of leaves is  totally disconnected; its complement is referred to as the skeleton.
One endows $\t$  with a probability measure $\mu$ carried by the subset of leaves and with a
sigma-finite length measure $\lambda$ carried by the skeleton. More precisely 
the distance between to points in $\t$, say $a,b$, is given by 
$\lambda([a,b])$, where $[a,b]$ stands for the path joining $a$ to $b$ in $\t$.

We next introduce marks on the skeleton of $\t$ that appear and disappear randomly as time passes at some constant rate. Specifically, we fix a parameter $r>0$
and introduce a Poisson point process on $\R\times (0,\infty)\times \t$ with intensity
$$\frac{1}{2} \d s \otimes \frac{1}{2r} \e^{-u/2r}\d u \otimes \d \lambda\,.$$
 An atom $(s,u,x)$ 
should be interpreted as follows: at time $s$ a mark appears at location $x$ on the skeleton and is erased at time $s+u$. In words, on any portion of the skeleton with length measure $\ell$,  marks appear with rate $\ell/2$ and the lifetime of each mark
is exponentially distributed with expectation $2r$, independently of the other marks.
For every $s\in\R$, we denote by
${\mathcal M}_{r,s}$ the random point measure on the skeleton of $\t$ induced by the marks present at time $s$. It is immediate to check that for each fixed $s$, conditionally on $\t$,  ${\mathcal M}_{r,s}$ is a Poisson point measure on the skeleton of $\t$ with intensity $r\lambda$. 

Following Aldous and Pitman \cite{AP}, for every $s\in\R$ we use the atoms of the Poisson random measure ${\mathcal M}_{r,s}$ to decompose the set of leaves of $\t$. More precisely, we decide that two leaves belong to the same component if and only if ${\mathcal M}_{r,s}$ has no mass on the path that joins those leaves
(note that the probability of this event is $\exp(-r\ell)$ where $\ell$ is the length
of the path between those two leaves). The components are closed in the subset of leaves, and we denote by
$\Y_{r,s}$ the sequence of their $\mu$-masses  ranked in the decreasing order. 
We stress that for every $s\in\R$, $\Y_{r,s}$ has the same distribution as $\F_r$, the CRT mass-fragmentation evaluated at time $r$, which was mentioned in the Introduction
and is described in Theorem 4 of \cite{AP}.

Recall from the Introduction the construction of the fragmentation-coagulation chain $\X\n_k$. We are now able to state

\begin{theorem}\label{T1} For each fixed $t>0$, the fragmentation-coagulation process 
$$(\X\n_{\lfloor tn + s\sqrt n\rfloor})_{s\in \R}$$
converges weakly in the sense of finite dimensional distributions as $n\to\infty$ towards 
the mixed process
$$(\Y_{R_t,s})_{s\in \R},$$
where $R_t$ denotes a random variable on $(0,\infty)$
which is independent of the preceding processes and has the distribution of a reflected Brownian motion at time $t$, viz.
$$\P(R_t\in \d r)= \sqrt {\frac{2}{t\pi}} \exp\left(-\frac{r^2}{2t}\right)\d r\,,\qquad r>0\,.$$
\end{theorem}

The rest of this section is devoted to the proof of Theorem \ref{T1}; the scheme of the argument is adapted from Aldous and Pitman \cite{AP}. We first recall the formulation of the convergence of uniform random trees towards the CTR via reduced trees.

 Given $\t$, we sample a sequence $U_1, \ldots$ of i.i.d. random leaves according to the law $\mu$, and for every integer $i\geq 1$, we denote by $\r(\infty,i)$ the subtree reduced to the first $i$ leaves, i.e. the smallest connected subset of $\t$ containing $U_1, \ldots, U_i$.  The reduced tree $\r(\infty,i)$ is a combinatorial tree (simple graph with no cycles) with leaves labeled by  $1,\ldots, i$  and some unlabeled internal nodes. The paths between two
adjacent internal nodes or between a leaf and an adjacent internal node are called edges. The
lengths of edges are given by the length measure $\lambda$ on $\t$, and the joint distribution of the shape and the edge-lengths is described by Lemma 21 in \cite{Al}.

 For every integer $n\geq 2$, we also consider a uniform random tree $\t_n$ on a set of $n$ vertices, say  $\{1, \ldots, n\}$, and assign length $1/\sqrt n$ to every edge. 
For every $1\leq i \leq n$, we denote by $\r(n,i)$ the sub-tree reduced to the first $i$ vertices. We agree that internal nodes with degree $2$ are discarded,
so that the edge-length between two adjacent vertices  in $\r(n,i)$ is $(1+j)/\sqrt n$ with $j$ the number of internal nodes with degree $2$ lying on the path connecting these vertices. 
It has been shown by Aldous  (see (49) in \cite{Al} or Lemma 9 in \cite{AP}) that 
for every fixed $i$, 
\begin{equation}\label{E1}
\r(n,i) \, \Longrightarrow \r(\infty,i)\quad\hbox{ as }n\to\infty\,,
\end{equation}
in the sense of weak convergence of the joint distributions of shape and edge lengths.

We then add and erase  marks randomly on $\t_n$ as explained in the Introduction.
For each integer $k\geq 0$, we denote by ${\mathcal M}\n_k$ a random point measure on $\t_n$ that assigns a unit mass to each marked edge after $k$ steps,
and for every $i\leq n$ by
${\mathcal M}^{(n,i)}_k$ the restriction of ${\mathcal M}\n_k$  to $\r(n,i)$.  Similarly, we also denote by ${\mathcal M}^{(\infty,i)}_{r,s}$
the restriction of the Poisson point measure ${\mathcal M}_{r,s}$  to the reduced tree $\r(\infty, i)$.
The key to Theorem \ref{T1} lies in the following limit theorem which can be viewed as a multi-dimensional extension of Equation (18) in \cite{AP}.
Essentially it is a
consequence of the law of rare events combined with the convergence of reduced trees.
Recall that the reduced trees are defined by their shapes and edge-lengths, 
and that the set of shapes of trees with $i$ vertices is finite. The reduced trees $\r(n,i)$ and $\r(\infty, i)$ should thus be viewed as random variables with values in some Polish space,
and weak convergence of random point measures should be understood in this setting.

\begin{lemma}\label{L1} Fix $t>0$ and an integer $i$. The process of random point measures on the reduced trees
$$({\mathcal M}^{(n,i)}_{\lfloor tn + s\sqrt n\rfloor}, \r(n,i))_{s\in \R}$$
converges weakly in the sense of finite dimensional distributions as $n\to\infty$ towards 
the mixed process of point measures on the reduced CRT
$$({\mathcal M}^{(\infty,i)}_{R_t,s}, \r(\infty,i))_{s\in \R},$$
where $R_t$ denotes a random variable on $(0,\infty)$
which is independent of the preceding processes and has the distribution of a reflected Brownian motion at time $t$, viz.
$$\P(R_t\in \d r)= \sqrt {\frac{2}{t\pi}} \exp\left(-\frac{r^2}{2t}\right)\d r\,,\qquad r>0\,.$$
\end{lemma}

\proof We first deal with the one-dimensional convergence in the statement,
rephrasing (and slightly developing) the argument for  Equation (18) in \cite{AP}.

By Skorohod's representation, we may assume that the convergence \eqref{E1} for the reduced trees holds almost surely and not merely in distribution. Thus with a probability close to $1$ when $n$ is large, the shape of $\r(n,i)$ coincides with that of $\r(\infty, i)$, and the edge lengths of $\r(n,i)$ and of $\r(\infty, i)$ are close.
We denote by $N\n_k$ the total number of marks on $\t_n$ after $k$ steps
and consider a sequence $(r_n)_{n\in\N}$ of integers with $r_n\sim r \sqrt n$ for some $r>0$. 
We first work for each $n$  conditionally on the event that 
$N\n_{\lfloor tn + s\sqrt n\rfloor}=r_n$.

 Recall that $\t_n$ has $n-1$ edges, each of length $1/\sqrt n$, so the number of edges 
 in a segment of $\t_n$ is $\sqrt n$ times the length of that segment. We stress that when $n$ is large, the number of edges in the reduced tree $\r(n,i)$ is of order $\sqrt n = o(n)$ and the number of marked edges in $\r(n,i)$ of order $O(1)= o(r_n)$. This is important to justify the claims of asymptotic independence which will be made below.
 
  As $(n-1)^{-1}r_n \sqrt n\sim r$, it follows from the law of rare events that when $n$ is large, the number of marked edges after $\lfloor tn + s\sqrt n\rfloor$ steps on a segment in $\t_n$ is
approximately Poisson distributed with parameter given by $r$ times the length of that segment, and further to disjoint segments correspond asymptotically independent Poisson variables. This entails that the conditional distribution of 
${\mathcal M}^{(n,i)}_{\lfloor tn + s\sqrt n\rfloor}$ converges weakly as $n\to\infty$
towards a Poisson random measure on $\r(\infty, i)$ with intensity $r \lambda$, i.e. 
\begin{equation}\label{E2}
{\mathcal L}\left({\mathcal M}^{(n,i)}_{\lfloor tn + s\sqrt n\rfloor}
\mid N\n_{\lfloor tn + s\sqrt n\rfloor}=r_n\right) \, \Longrightarrow \, {\mathcal M}^{(\infty,i)}_{r,s}\,,
\end{equation}
where the notation ${\mathcal L}(Z\mid \Lambda)$ refers to the conditional law of the variable $Z$ given the event $\Lambda$.

We next present the main lines of the argument for extending \eqref{E2} to multi-dimensional convergence by analyzing the evolution of the random point measures as $s$ increases. 
It is readily checked  that with probability one
$$N\n_{\lfloor tn + (s+s')\sqrt n\rfloor} \sim r_n \sim r\sqrt n $$
uniformly for $s'\geq 0$ in an arbitrary bounded interval. Thus for every
$k=nt + O(\sqrt n)$, each atom of ${\mathcal M}^{(n,i)}_k$ has a probability
close to $1/(2 r_n)\sim 1/(2r\sqrt n)$ to be erased at the next step, where the factor $1/2$ accounts for the probability that head turns up when the fair coin is flipped. 
The probability that a given atom of 
 ${\mathcal M}^{(n,i)}_{\lfloor tn + s\sqrt n\rfloor}$  has not been erased after $\lfloor s'\sqrt n\rfloor$ further steps is close to
$$\left(1-1/(2r_n)\right)^{\lfloor s'\sqrt n\rfloor}\sim \exp(-s'/(2r))$$
when $n$ is large; in other words if one unit of time corresponds to $\sqrt n$ steps, each atom of ${\mathcal M}^{(n,i)}_{\lfloor tn + s\sqrt n\rfloor}$ is removed after a time which is approximately exponentially distributed with mean $2r$. 
A similar argument shows that asymptotically, each atom of  ${\mathcal M}^{(n,i)}_{\lfloor tn + s\sqrt n\rfloor}$ is removed or not after $\lfloor s'\sqrt n\rfloor$ more steps
 independently on the other atoms.

On the other hand, at each step, a mark appears on an un-marked edge with probability 
close to $1/(2(n-r_n))\sim 1/(2n)$. Recalling that any given mark is also erased at each step with probability close to $1/(2r\sqrt n)$, and neglecting the event of multiple appearances and erasures of a mark whose probability is of lower order, 
we deduce that the probability that an edge with no mark after ${\lfloor tn + s\sqrt n\rfloor}$ steps be marked after ${\lfloor s'\sqrt n\rfloor}$ further steps is 
\begin{eqnarray*}
\sum_{j=1}^{{\lfloor s'\sqrt n\rfloor}} \frac{1}{2n}\left(1-\frac{1}{2r\sqrt n}\right)^j
&= &\frac{r}{\sqrt n}\left( \left(1-\frac{1}{2r\sqrt n}\right)- \left(1-\frac{1}{2r\sqrt n}\right)^{\lfloor s'\sqrt n\rfloor +1}\right) \\
&\sim& \frac{r}{\sqrt n}\left(
1-\exp\left(-\frac{s'}{2r}\right) \right)\,.
\end{eqnarray*}
It then follows from the law of rare events that for
any given segment of $\r(n,i)$ with length $\ell$ (i.e. with $\ell \sqrt n$ edges), the 
number of marked edges after  ${\lfloor tn + (s+s')\sqrt n\rfloor}$  steps which were un-marked after  ${\lfloor tn + s\sqrt n\rfloor}$  steps is approximately Poisson
with parameter 
$$r\ell\left(
1-\exp\left(-\frac{s'}{2r}\right) \right)\,. $$
Further, one checks readily that the evolutions of marks on a given finite sequence of disjoint segments are asymptotically independent.   

Putting the pieces together, this shows that when $n$ is large, the distribution of 
${\mathcal M}^{(n,i)}_{\lfloor tn + (s+s')\sqrt n \rfloor}$ given $ {\mathcal M}^{(n,i)}_{\lfloor tn + s\sqrt n\rfloor}$ is close to that of a measure obtained from $ {\mathcal M}^{(n,i)}_{\lfloor tn + s\sqrt n\rfloor}$ by removing each atom with probability $1-\exp(-s'/(2r))$
independently one of the others (i.e. by thinning), and further adding an independent 
Poisson measure on $\r(n,i)$ with intensity  $r\left(
1-\exp\left(-s'/{2r}\right) \right)\lambda$. 
Comparing with the evolution of the random point measure ${\mathcal M}^{(\infty,i)}_{r,s}$ when $s$ increases, we see that  \eqref{E2} can be extended as follows: 
for every $s'\geq 0$ we have
$$
{\mathcal L}\left(\left( {\mathcal M}^{(n,i)}_{\lfloor tn + s\sqrt n\rfloor}, {\mathcal M}^{(n,i)}_{\lfloor tn + (s+s')\sqrt n\rfloor}\right)
\mid N\n_{\lfloor tn + s\sqrt n\rfloor}=r_n\right) \, \Longrightarrow \,
\left( {\mathcal M}^{(\infty,i)}_{r,s}, {\mathcal M}^{(\infty,i)}_{r,s+s'}\right)\,.
$$
An iteration based on the obvious Markov property of the processes of the radom measures involved (in the case of $({\mathcal M}^{(\infty,i)}_{r,s})_{s\in\R}$, this Markov property follows from the absence of memory of the exponential variables that specify the lifetime of the marks) enables us to conclude that 
\begin{equation}\label{E3}
{\mathcal L}\left(\left({\mathcal M}^{(n,i)}_{\lfloor tn + s'\sqrt n\rfloor}\right)_{s'\geq s}
\mid N\n_{\lfloor tn + s\sqrt n\rfloor}=r_n\right) \, \Longrightarrow \,
\left( {\mathcal M}^{(\infty,i)}_{r,s'}\right)_{s'\geq s}\,,
\end{equation}
in the sense of finite dimensional distributions.

All that we need now is to get rid of the conditioning in \eqref{E3}, which is straightforward. Indeed Donsker's invariance principle shows that there is the weak convergence $$\frac{1}{\sqrt n}N\n_{\lfloor tn + s\sqrt n\rfloor}\Longrightarrow R_t\,,$$
where  $R_t$ is a random variable distributed as in the statement. 
On the other hand, it is easy seen from the construction of the random point measures ${\mathcal M}_{r,s}$
that the finite-dimensional distributions of the process $ \left( {\mathcal M}^{(\infty,i)}_{r,s'}\right)_{s'\geq s}$ depend continuously on the parameter $r$.
We derive from above that 
$$\left({\mathcal M}^{(n,i)}_{\lfloor tn + s'\sqrt n\rfloor}\right)_{s'\geq s}
 \, \Longrightarrow \,
\left( {\mathcal M}^{(\infty,i)}_{R_t,s'}\right)_{s'\geq s}$$
in the sense of finite dimensional distributions, which is our claim. \QED

We now explain how Theorem \ref{T1} follows from Lemma \ref{L1}, focussing on
one-dimensional distributions as the multi-dimensional case is similar but with heavier notation.

For every integers $n,k\geq 1$, we denote by $\Pi\n_k$ the random partition of $\{1,\ldots, n\}$
induced by the marks on edges of $\t_n$ after $k$ steps, that is the blocks of $\Pi\n_k$
are characterized by the property that there is no mark on the paths connecting two vertices 
in the same block.  By definition, $\X\n_k$ is the sequence of the sizes of the blocks of $\Pi\n_k$ ranked in the decreasing order and rescaled by a factor $1/n$.

Similarly, 
for every  $s\in\R$, we denote by
$\Pi^{(\infty)}_{r,s}$ the random partition of $\N$ such that two integers $j,j'$ belong to the same block of $\Pi^{(\infty)}_{r,s}$ if and only if  there is no atom of
 the independent  mixed-Poisson random measure ${\mathcal M}_{r,s}$ on the path in $\t$ from the leaf $U_j$ to the leaf $U_{j'}$. Plainly the random partition $\Pi^{(\infty)}_{r,s}$ is exchangeable and the asymptotic frequencies of its blocks are
given by $\Y_{r,s}$, i.e. the  sequence of 
the $\mu$-masses of the connected components of the CRT $\t$ cut at the atoms of ${\mathcal M}_{r,s}$.

For every integer $i\leq n$, we also denote by $\Pi^{(n,i)}_k$ (respectively by  $\Pi^{(\infty,i)}_{r,s}$)
the restriction of $\Pi^{(n)}_k$   (respectively of   $\Pi^{(\infty)}_{r,s}$)  
to the first $i$ vertices. Plainly these restricted partition only depend on the reduced tree $\r(n,i)$ and $\r(\infty,i)$, and the marks on their edges after $k$ steps and the atoms of
the random measure ${\mathcal M}^{(\infty,i)}_{r,s}$, respectively. Lemma \ref{L1} implies that for every $i$, in the obvious notation, when $n\to \infty$
 there is the weak convergence
\begin{equation}\label{E4}
\Pi^{(n,i)}_{\lfloor tn + s \sqrt n\rfloor} \, \Longrightarrow\, \Pi^{(\infty,i)}_{R_t,s}\,.
\end{equation}

Repeating the  argument of Aldous and Pitman for proving Theorem 3 in \cite{AP}
enables us to conclude from \eqref{E4} and the preceding observations that
$$\X\n_{\lfloor tn + s \sqrt n\rfloor} \, \Longrightarrow\, \Y_{R_t,s}\,,$$
which is the one-dimensional version of Theorem \ref{T1}. The multidimensional case is similar, using the full strength of Lemma \ref{L1}.

\end{section}

\begin{section}{A comment, a complement, and an open question}

The stationary limiting process $(\Y_{R_t,s})_{s\in\R}$ which appears in Theorem \ref{T1} is expressed as a mixture. The mixing variable $R_t$ may be thought of as the effective age of the system as it represents the intensity of cuts along the skeleton of the CRT. In this direction, we note that the variables $R_t$ are stochastically increasing with $t$. We also mention that $R_t$ can be recovered from a sample of 
$\Y_{R_t,s}$. Indeed, it follows easily from Theorem 4 in Aldous and Pitman \cite{AP}
and the law of large numbers for Poisson processes that if $F_t(i)$ denotes the
$i$-th largest term of $\F_t$, the CRT fragmentation of masses observed at time $t$, then
with probability one
$$F_t(i)\sim \frac{2}{\pi t^2 i^2}\qquad\hbox{ as }i\to \infty\,.$$
Equivalently 
$$Y_{r,s}(i)\sim \frac{2}{\pi r^2 i^2}\qquad\hbox{ as }i\to \infty\,,$$
where $Y_{r,s}(i)$ denotes the $i$-th largest term of $\Y_{r,s}$
and we conclude that
$$R_t=\lim_{i\to\infty} \sqrt{\frac{2}{\pi Y_{R_t,s}(i) i^2}}
\qquad \hbox{a.s.}$$

We used uniform random trees $\t_n$ merely to stick to Pitman's original framework \cite{Pi}. Nonetheless the same results hold if we replace $\t_n$ by any other sequence of discrete random tree which converges to the Brownian CRT after rescaling edge lengths by a factor $1/\sqrt n$ and masses of vertices by a factor $1/n$ (for instance critical Galton-Watson 
trees with finite variance and conditioned to have total size $n$).
In the same vein, the results of this note can be extended to certain sequences of so-called birthday trees. Indeed, Camarri and Pitman \cite{CP} have established the weak convergence of suitably rescaled birthday trees towards certain Inhomogeneous Continuum Random Trees. On the other hand, dynamics of edge-deletion for birthday trees
bears the same connection to the additive coalescence as uniform random trees, except that the initial distribution of masses is inhomogeneous. The asymptotic behavior of the latter has been characterized by Aldous and Pitman \cite{AP2}, in the study of the entrance boundary of the additive coalescence. We thus have all the ingredients needed to apply the arguments of the present work to this more general setting. Of course, the limiting processes will then have different distributions. 

Our aim in this work was to point at the phenomenon of subaging in a fragmentation-coagulation process. The model that we used for this purpose is easy to deal with
although somewhat artificial. There are other  discrete models for the evolution of random forests which may be more natural, but are also much harder to investigate. Here is an example, which inspired by the subtree prune and regraft algorithm; see Chapter 9 of Evans \cite{Ev}. We now work with rooted forests on $n$ vertices, that is each tree has one distinguished vertex that serves as the root. At each step we flip a fair coin.
With probability $1/2$ we delete an edge chosen uniformly at random in this forest. This disconnects the tree containing that edge into two rooted subtrees. With probability $1/2$, we create a new edge between a vertex chosen uniformly at random and the root of a tree chosen uniformly at random amongst the trees to which the chosen vertex does not belong. Our result suggests that a similar subaging phenomenon might occur at the same scale as in the present study. Proving or disproving this property would be interesting, but does not seem easy.

\end{section}

\noindent{\bf Acknowledgments}.  I would like to thank an Associate Editor and two anonymous referees for their pertinent comments and suggestions on the first draft. This work has been supported by ANR-08-BLAN-0220-01.

  \end{document}